\documentclass[11pt]{article}

\usepackage{amsfonts, amsmath, amsthm, amssymb}
\usepackage{verbatim}
\usepackage{graphicx}

\date{}
\title{\vspace{-0.8cm}On the size of minimal unsatisfiable formulas \thanks{This research forms part of the Ph.D thesis written by the author under the supervision of Prof. Benny Sudakov.}}
\author{
Choongbum Lee \thanks{Department of Mathematics, UCLA, Los Angeles, CA, 90095. E-mail:
choongbum.lee@gmail.com.
Research supported in part by Samsung Scholarship.
}
}

\newtheorem{thm}{Theorem}
\newtheorem{prop}[thm]{Proposition}
\newtheorem{lemma}[thm]{Lemma}

\newtheorem{claim}[thm]{Claim}

\newtheorem*{ques}{Question}

\begin{document}
\maketitle

\begin{abstract}
An unsatisfiable formula is called minimal if it becomes satisfiable whenever any of its clauses are removed. We construct minimal unsatisfiable $k$-SAT formulas with $\Omega(n^k)$ clauses for $k \geq 3$, thereby negatively answering a question of Rosenfeld. This should be compared to the result of Lov\'asz \cite{Lovasz1} which asserts that a critically 3-chromatic $k$-uniform hypergraph can have at most $\binom{n}{k-1}$ edges. \\
\end{abstract}

\section{Introduction}

Given $n$ boolean variables $x_1, \ldots, x_n$, a \textit{literal} is a variable $x_i$ or its negation $\overline{x}_i \, (1 \leq i \leq n)$. A \textit{clause} is a disjuction of literals and by \textit{$k$-clause} we denote a clause of size $k$. A \textit{CNF(Conjunctive Normal Form) formula} is a conjunction of clauses and a \textit{$k$-SAT formula} is a CNF formula with only $k$-clauses. Throughout this article \textit{formula} will mean a CNF formula and it will be given as a pair $F=(V,C)$ with variables $V=\{x_1, \ldots, x_n\}$ and  clauses $C$ as collection of disjunction of literals $V \cup \overline{V}$. A formula is called \textit{satisfiable} if there exists an assignment of values to variables so that the formula becomes true. A formula is called \textit{minimal unsatisfiable} if it is not satisfiable but removing any clause makes it satisfiable. 

Satisfiablity of a formula is closely related to the 2-colorability of a hypergraph in the following sense. A formula is satisfiable if there is an assignment of values to variables in a way that no clauses have only false literals inside it. Similarily a hypergraph is 2-colorable if there is a way to color the vertices into two colors so that none of the edges become monochromatic. A hypergraph $H=(V,E)$ is called critically 3-chromatic if it is not 2 colorable but the deletion of any edge makes it 2 colorable. In this analogy, minimal unsatisfiable formulas correspond to critically 3-chromatic hypergraphs. Therefore it is natural to ask if similar results hold for both problems. In particular, we are interested if the same restriction on the number of clauses (edges, respectively) holds or not.

In the case of lower bounds, roughly the same estimate holds for both formulas and hypergraphs. Seymour \cite{Seymour1} used linear algebra method to deduce that a critically 3-chromatic hypergraph $H=(V,E)$ must satisfy $|E| \geq |V|$ if there is no isolated vertex. The corresponding bound for CNF formulas appeared in Aharoni and Linial \cite{AharoniLinial1} where they proved that minimal unsatisfiable formula $F=(V,C)$ must satisfy $|C| \geq |V|+1$ if every variable is contained in some clause.

For uniform hypergraphs there are also known upper bound results.  Lov\'asz \cite{Lovasz1} proved that any critically 3-chromatic $k$-uniform hypergraph has at most $\binom{n}{k-1}$ edges. This result is asymptotically tight, as was shown by Toft \cite{Toft1} who constructed critically 3-chromatic $k$-uniform hypergraphs with $\Omega(n^{k-1})$ edges.

Motivated by these results, Rosenfeld \cite{Rosenfeld1} asked if the analogy also holds for minimal unsatisfiable $k$-SAT formulas. 

\begin{ques}
Should minimal unsatisfiable $k$-SAT formulas have at most $O(n^{k-1})$ clauses?
\end{ques}

It is not difficult to show that this conjecture is true for $k=2$ and we will give the simple proof of this in section \ref{2SATsection}. However for $k$-SAT formulas with $k \geq 3$ we show that surprisingly the answer for the question of Rosenfeld is negative. In section 3 we will construct minimal unsatisfiable $k$-SAT formulas with $\Omega(n^k)$ clauses.

\section{2-SAT formulas}
\label{2SATsection}

First we give explicit minimal unsatisfiable 2-SAT formulas. Consider the 2-SAT formula $F^{(2)} = (V^{(2)}, C^{(2)})$ where $V^{(2)}=\{ y_1, y_2, \ldots, y_{2l}\}$ and $C^{(2)}= \{ y_i \vee y_{i+1}, \overline{y}_i \vee \overline{y}_{i+1} : i=1,2,\ldots, 2l-1 \} \cup \{y_1 \vee \overline{y}_{2l} \} \cup \{\overline{y}_1 \vee y_{2l}\}$. $F^{(2)}$ is unsatisfiable because if $y_i = y_{i+1}$ for some $i$ then either $y_i \vee y_{i+1}$ or $\overline{y}_i \vee \overline{y}_{i+1}$ is false and otherwise if $y_i = \overline{y}_{i+1}$ for all $1 \leq i \leq 2l-1$ then $y_1 = \overline{y}_{2l}$ and this time either $y_1 \vee \overline{y_{2l}}$ or $\overline{y_1} \vee y_{2l}$ will become false.
 To prove that $F^{(2)}$ is minimal unsatisfiable, we only check that deleting $y_1 \vee y_2$ or $y_1 \vee \overline{y}_{2l}$ makes the new formula satisfiable as other clauses can be checked similarily. 
 In each case, the assignment of ($y_1 = y_2 = \textrm{false}$, $y_3 = \ldots = y_{2l-1} = \textrm{true}$, $y_4 = \ldots = y_{2l} = \textrm{false}$) and ($y_1 = y_3 = \ldots = y_{2l-1} = \textrm{false}$, $y_2 = y_4 = \ldots = y_{2l} = \textrm{true}$) will make the remaining clauses true.

Next we prove the linear upper bound of number of clauses in minimal unsatisfiable 2-SAT formulas.
\begin{prop} 
\label{2SATprop}
Minimal unsatisfiable 2-SAT formulas have at most $4n$ clauses.
\end{prop}

\begin{proof}

Given a minimal unsatisfiable 2-SAT formula $F = (V, C)$, let's consider the \textit{implication graph} $D$ of this 2-SAT formula which is the directed graph $D$ over the vertices $V \cup \overline{V}$ with two directed edges corresponding to each clause $z_1 \vee z_2 \in C$  given as $\overline{z}_1 \rightarrow z_2$ and $\overline{z}_2 \rightarrow z_1$. Aspvall, Plass and Tarjan \cite{AspvallPlassTarjan1} proved that 2-SAT is unsatisfiable if and only if its implication graph has a strongly connected component which contains both $x_i$ and $\overline{x}_i$ for some index $i$. Therefore the unsatisfiability of $F$ implies the existence of directed path from $x_i$ to $\overline{x}_i$ and from $\overline{x}_i$ and $x_i$ in $D$ for some index $i$. Now observe that the minimality of $F$ forces every clause $z_1 \vee z_2 \in C$ to have at least one of its corresponding edge $\overline{z}_1 \rightarrow z_2$ or $\overline{z}_2 \rightarrow z_1$ in these directed paths. As otherwise deleting the clause will not change the unsatisfiability of $F$ (because it still contains both directed paths). Since $D$ has $2n$ vertices, there can be at most $4n$ edges in the two directed paths. Therefore $|C| \leq 4n$.
\end{proof}

\section{$k$-SAT formulas}
\label{3SATsection}

In this section we construct minimal unsatisfiable $k$-SAT formulas on $n$ variables with $\Omega(n^k)$ clauses. For simplicity we describe in details the construction of 3-SAT formulas only. This construction can be easily generalized for all $k$. Informally, start with a minimal unsatisfiable ``almost'' 3-SAT formula with $\Omega(n^3)$ clauses where ``almost'' means that only a small number of clauses is not of size 3. Then transform this formula into a ``genuine'' 3-SAT formula by replacing the clauses of size greater than 3 by 3-clauses while keeping the minimal unsatisfiable property. During the process the number of variables will not increase too much and therefore we will end up with a 3-SAT formula that we have promised. Now we should make it into a formal argument.

The following lemma will allow us to change the size of a clause in the formula. This lemma is a modified version of Theorem 1 and 4 in \cite{Toft1} which were originally used by Toft to construct $k$-uniform hypergraphs with $\Omega(n^{k-1})$ edges.

Let $F_X = (V_X, C_X), F_Y = (V_Y, C_Y)$ be formulas with disjoint sets of variables(that is, $V_X \cap V_Y = \emptyset$) and $c_0 = z_1 \vee z_2 \vee \ldots z_k \in C_X$ be a $k$-clause of $F_X$ where $k \leq |C_Y|$. For an arbitrary surjective map $h$ from $C_Y$ to $\{z_1, z_2, \ldots, z_k \}$, let the formula $F_Z = (V_Z, C_Z)$ be as following.
\begin{itemize}
\item $V_Z = V_X \cup V_Y$ \qquad $C_Z = (C_X \backslash \{c_0\}) \cup \{ c_y \vee h(c_y)  : c_y \in C_y \} $
\end{itemize}

\begin{lemma}
\label{replacelemma}
If $F_X$ and $F_Y$ are minimal unsatisfiable formulas, then $F_Z$ constructed as above is also a minimal unsatisfiable formula.
\end{lemma}
\begin{proof}
Let's first show that $F_Z$ is unsatisfiable. For arbitrary values of $V_X$ there must exist a clause $c_x \in C_X$ which is false. If $c_x \neq c$ then we are done as $c_x \in C_Z$ so assume that $c_x = c$. Since every literal $x \in c_0$ is false, a new clause $c_y \vee h(c_y)$ is true if and only if $c_y$ is true. But $F_Y$ is unsatisfiable so there must exist a clause $c_y$ which is false and therefore $F_Z$ is unsatisfiable. 

Next we prove that removing any clause $c_z \in C_Z$ makes $F_Z$ satisfiable. First asuume that $c_z \in C_X \backslash \{c_0\}$. Then give values to $V_X$ so that every clause in $C_X$ except $c_z$ is satisfied. Since $c_z \neq c_0$, there must exist a literal $x \in c_0$ which is true. Pick a clause $c' \in h^{-1}(x) \subset C_Y$ ($h^{-1}(x)$ is non-empty because $h$ is surjective) and give $V_Y$ the values which make every clause except $c'$ in $C_Y$ true. Observe that every clause in $C_X \backslash \{c_0\}$ except $c_z$ is true by values of $V_X$ and every clause in $\{ c_y \vee h(c_y) : c_y \in C_y \}$ is true either by values of $V_Y$ or the literal $x$. Now assume that $c_z = c' \vee x \in \big\{ c_y \vee h(c_y)  : c_y \in C_y \big\}$ and give $V_X$ the values which make every clause except $c_0$ true and give $V_Y$ the values which makes every clause except $c'$ true. This assignment of values will make every clause but $c_z \in C_Z$ true and thus we are done.
\end{proof}

Next step is to construct an ``almost'' 3-SAT formula with many clauses. Let $V_0 = \{ x_1, x_2, \ldots, x_{6m} \}$ and look at the formula $F_0 = (V_0, C_0)$ with clauses given as, 

\begin{itemize}
\item $C_0 = \{ x_{i_1} \vee x_{i_2} \vee x_{i_3} : 1 \leq i_1 \leq 2m, 2m+1 \leq i_2 \leq 4m, 4m+1 \leq i_3 \leq 6m\}$ \\ 
\hspace*{0.8cm} $\cup \{ \overline{x}_1 \vee \overline{x}_2 \vee \ldots, \overline{x}_{2m} \} \cup \{ \overline{x}_{2m+1}\vee \ldots \vee \overline{x}_{4m} \} \cup \{ \overline{x}_{4m+1} \vee \ldots \vee \overline{x}_{6m} \}$
\end{itemize}

Informally, partition the variables $V$ into three equal parts $V_1, V_2, V_3$ and consider every clauses $x_1 \vee x_2 \vee x_3$ with $x_i \in V_i$ and add three more clauses $\overline{V}_1, \overline{V}_2, \overline{V}_3$. Note that this formula contains $(2m)^3 + 3$ clauses.

\begin{claim} 
$F_0$ is a minimal unsatisfiable formula.
\end{claim}
\begin{proof}
Let's first prove that $F_0$ is unsatisfiable. Assume that the three clauses $\overline{x}_1 \vee \overline{x}_2 \vee \ldots \vee \overline{x}_{2m}$, $\overline{x}_{2m+1} \vee \ldots \vee \overline{x}_{4m}$, $\overline{x}_{4m+1} \vee \ldots \vee \overline{x}_{6m} $ are all true. Then there must exist $1 \leq i_1 \leq 2m, 2m+1 \leq i_2 \leq 4m, 4m+1 \leq i_3 \leq 6m$ such that $x_{i_1}=x_{i_2}=x_{i_3}=\textrm{false}$. But this will make the clause $x_{i_1} \vee x_{i_2} \vee x_{i_3}$ false. Therefore $F_0$ is unsatisfiable.

Now assume that we remove a clause $c$. If $c = x_{i_1} \vee x_{i_2} \vee x_{i_3}$ for some $i_1, i_2, i_3$ then assigning $x_{i_1}=x_{i_2}=x_{i_3}=\textrm{false}$ and everything else true will make the remaining part satisfiable. On the other hand if $c = \overline{x}_1 \vee \overline{x}_2 \vee \ldots \vee \overline{x}_{2m}$ then assigning $x_1=x_2=\ldots=x_{2m}=\textrm{true}$ and everything else false will make the remaining part satisfiable. Similar assignment will work for clauses $\overline{x}_{2m+1} \vee \ldots \vee \overline{x}_{4m}$ and $\overline{x}_{4m+1} \vee \ldots \vee \overline{x}_{6m}$.
\end{proof}

\noindent \textbf{Construction}

Note that the formula $F_0$ is ``almost'' a 3-SAT formula in the sense that there are only three clauses whose size is not 3. Use Lemma \ref{replacelemma} with $F_X=F_0$, $c_0=\overline{x}_1 \vee \overline{x}_2 \vee \ldots \vee \overline{x}_{2m} \in C_0$ and $F_Y=F^{(2)}$ where $F^{(2)}$ is a minimal unsatisfiable 2-SAT formula with $m$ variables and $2m$ clauses as constructed in section \ref{2SATsection}. The obtained formula $F_1$ is a minimal unsatisfiable formula over $6m+m = 7m$ variables and has only two clauses whose size are not 3. (All new clauses are 3-clauses.)

Repeat the same process with the remaining two $2m$-clauses to obtain a minimal unsatisfiable formula $F_2$ whose every clause has size 3 i.e. $F_2$ is a 3-SAT formula over $n=9m$ variables. Note that it still contains the original 3-clauses $\{ x_{i_1} \vee x_{i_2} \vee x_{i_3} : 1 \leq i_1 \leq 2m, 2m+1 \leq i_2 \leq 4m, 4m+1 \leq i_3 \leq 6m\}$. There are $8m^3 = (\frac{2}{9}n)^3$ such clauses and therefore this 3-SAT formula $F_2$ contains $\Omega(n^3)$ clauses.

For $k \geq 4$, minimal unsatisfiable $k$-SAT formulas with $\Omega(n^k)$ clauses can be constructed similarily. Use $F_0^{(k)} = (V_0^{(k)}, C_0^{(k)})$ where,
\begin{itemize}
\item $V_0^{(k)} = \{ x_1, x_2, \ldots, x_{mk} \}$
\item $C_0^{(k)} = \{ x_{i_1} \vee x_{i_2} \vee \ldots \vee x_{i_k} \} : (t-1)m + 1 \leq i_t \leq tm, 1 \leq t \leq k \}$ \\
\hspace*{0.9cm} $\cup \big(\cup_{s=1}^{k} \{ \overline{x}_{(s-1)m + 1} \vee \overline{x}_{(s-1)m + 2} \vee \ldots \vee \overline{x}_{sm} \}\big)$
\end{itemize}
By the same process as above one can verify that $F_0^{(k)}$ is minimal unsatisfiable. Then replace the $m$-clauses by $k$-clauses using Lemma \ref{replacelemma} and minimal unsatisfiable $(k-1)$-SAT formulas. The final formula will be a minimal unsatisfiable $k$-SAT formula with $\Omega(n^k)$ clauses. Details are omitted. \\

\noindent \textbf{Acknowledgement.} I gratefully thank Benny Sudakov for his advice and guidance. I am also thankful to Po-Shen Loh for his useful propositions and corrections and to Boris Bukh for the fruitful discussion.

\end{document}